\documentclass[12pt]{amsart}
\usepackage{amsmath}
\usepackage{amsthm}
\usepackage{amssymb}
\usepackage{tikz}
\usepackage{hyperref}

\usepackage{url}
\usepackage{cite}
\usepackage{dsfont}
\usepackage{array}

\usepackage[all]{xy}
\usepackage{mathrsfs}
\usepackage{amscd}
\usepackage{latexsym}
\usepackage{epsfig}
\usepackage{graphicx}
\usepackage{amsfonts}
\usepackage{psfrag}
\usepackage{caption}
\usepackage{rotating}
\usepackage{mathtools}
\usepackage{etoolbox}

\usepackage{bm}

\newtheorem{theorem}{Theorem}
\newtheorem{proposition}[theorem]{Proposition}

\theoremstyle{definition}
\newtheorem{definition}[theorem]{Definition}
\newtheorem{assumption}[theorem]{Assumption}

\theoremstyle{remark}
\newtheorem*{remark}{Remark}

\DeclareMathOperator{\Ad}{Ad}
\DeclareMathOperator{\ad}{ad}
\DeclareMathOperator{\Cd}{Cd}
\DeclareMathOperator{\cd}{cd}

\DeclareMathOperator{\SU}{SU}
\DeclareMathOperator{\End}{End}
\DeclareMathOperator{\Aut}{Aut}
\DeclareMathOperator{\Stab}{Stab}
\DeclareMathOperator{\pr}{pr}

\DeclareMathOperator{\sgn}{sgn}

\newcommand{\orb}{\mathcal{O}}
\newcommand{\lieg}{\mathfrak{g}}
\newcommand{\liegdual}{\mathfrak{g}^{*}}
\newcommand{\liet}{\mathfrak{t}}
\newcommand{\lietdual}{\mathfrak{t}^{*}}

\newcommand{\red}{\mathrm{red}}

\newcommand{\CC}{\bf C}
\newcommand{\RR}{\bf R}
\newcommand{\ZZ}{\bf Z}

\begin{document}

\title[Intersection Pairings in $N$-fold Reduced Products]{Intersection Pairings in the  $N$-fold Reduced Product of Adjoint Orbits}

\author{Lisa C. Jeffrey}
\address{Department of Mathematics \\
University of Toronto \\ Toronto, Ontario \\ Canada}
\email{jeffrey@math.toronto.edu}
\urladdr{\url{http://www.math.toronto.edu/~jeffrey}} 

\author{Jia Ji}
\address{Department of Mathematics \\
University of Toronto \\ Toronto, Ontario \\ Canada}
\email{jia.ji@mail.utoronto.ca}

\thanks{The first author is partially supported by an NSERC Discovery Grants. The authors wish to thank Rebecca Goldin and Augustin-Liviu Mare for helpful conversations.}


\date{\today}


\begin{abstract}
In previous work we computed the symplectic 
volume of the symplectic  reduced space of the product 
of $N$ adjoint orbits of a compact Lie group. In this paper
we compute the intersection pairings in the same object.
\end{abstract}

\maketitle

 
\makeatletter
\def\@tocline#1#2#3#4#5#6#7{\relax
  \ifnum #1>\c@tocdepth 
  \else
    \par \addpenalty\@secpenalty\addvspace{#2}%
    \begingroup \hyphenpenalty\@M
    \@ifempty{#4}{%
      \@tempdima\csname r@tocindent\number#1\endcsname\relax
    }{%
      \@tempdima#4\relax
    }%
    \parindent\z@ \leftskip#3\relax \advance\leftskip\@tempdima\relax
    \rightskip\@pnumwidth plus4em \parfillskip-\@pnumwidth
    #5\leavevmode\hskip-\@tempdima
      \ifcase #1
       \or\or \hskip 1em \or \hskip 2em \else \hskip 3em \fi%
      #6\nobreak\relax
    \hfill\hbox to\@pnumwidth{\@tocpagenum{#7}}\par
    \nobreak
    \endgroup
  \fi}
\makeatother

\setcounter{tocdepth}{5}
\tableofcontents

\section{Introduction} \label{s:1}
Let $G$ be  a compact connected Lie group  with 
maximal torus $T$.
As a vector space, 
the equivariant cohomology of  a Hamiltonian $G$-space $M$ 
     is isomorphic to  the tensor product of the ordinary 
cohomology of $M$   and the $G$-equivariant cohomology of a point.
%
Here $S(\liet)$ is the polynomial ring on the Lie algebra of the maximal
torus $T$, which is denoted $\liet$.
This result comes from \cite{Kthesis} (Proposition  5.8). 
The above isomorphism is 
only an isomorphism of vector spaces, not of rings.

When $M$ and $G$ are as above,
 there is a surjective ring homomorphism $\kappa$ (the Kirwan map)
 from the equivariant
cohomology of $M$ to the ordinary cohomology 
of the symplectic reduced space or 
symplectic quotient 
$M_{\rm red} $, which is defined as 
$$M_{\rm red} = \mu^{-1}(0)/G$$
 where $\mu$ is the 
moment map. The ordinary cohomology of the reduced space
is the quotient of the equivariant cohomology of $M$ by the kernel of $\kappa$.

Provided  the reduced space is a smooth manifold, it satisfies
Poincar\'e duality, so its cohomology ring is determined by 
the intersection pairings (in other words the evaluation of 
cohomology classes against the fundamental class).

Let $M $ be the product of a collection of adjoint orbits of $G$.
In this situation, the above isomorphism is an isomorphism of 
$H^*_G({\rm pt})$-modules.
We give a formula  for the intersection 
pairings  in $M_{\rm red}$ 
using the same methods as 
in our earlier paper \cite{JJ18a}, in other words the localization 
theorem of Atiyah-Bott and Berline-Vergne and 
the residue formula of \cite{JK95} (Theorem 8.1).


%
%

\section{Notation and Conventions}\label{SC001}

Let $G$ be a compact connected Lie group. Let $\lieg$ be the Lie algebra of~ $G$. Let $\liegdual$ be the dual vector space of $\lieg$.

We choose a maximal torus $T$ in $G$. Let $\liet$ be the Lie algebra of~ $T$. Let $\lietdual$ be the dual vector space of $\liet$. Let $W = N_{G}(T)/T$ be the corresponding Weyl group.

Let $\Ad: G \to \Aut(\lieg)$ be the adjoint representation of $G$. Let $\Cd: G \to \Aut(\liegdual)$ be the coadjoint representation of $G$. More explicitly,
\begin{equation}
\left\langle \Cd(g) \xi, X \right\rangle = \left\langle \xi, \Ad(g^{-1}) X \right\rangle
\end{equation}
for all $g \in G$, $X \in \lieg$, $\xi \in \liegdual$, where $\left\langle \cdot, \cdot \right\rangle$ is the natural pairing between a covector and a vector.

\begin{remark}
Note that for all $g, h \in G$, $\Ad(g) \circ \Ad(h) = \Ad(gh)$ and $\Cd(g) \circ \Cd(h) = \Cd(gh)$. 
That is, both $\Ad$ and $\Cd$ are left actions. 
\end{remark}

Let $\ad: \lieg \to \End(\lieg)$ be the adjoint representation of $\lieg$. Let $\cd: \lieg \to \End(\liegdual)$ denote the coadjoint representation of the Lie algebra $\lieg$. Thus, $\cd(X) = - \ad(X)^{*}$.

\begin{remark}
Note that both $\ad: \lieg \to \End(\lieg)$ and $\cd: \lieg \to \End(\liegdual)$ are Lie algebra homomorphisms.
\end{remark}

For convenience we work with orbits of the adjoint action rather than the coadjoint action, so our orbits are subsets of $\lieg$ instead of $\liegdual$.
The invariant inner product $< \cdot, \cdot > $  on $\lieg$ (invariant under the adjoint action) gives a 
$G$-equivariant isomorphism  between $\lieg$  (equipped with the adjoint action) and $\liegdual$ (with the coadjoint action).


Let $\orb(\xi)$ denote the adjoint orbit through $\xi \in \lieg$. The following theorem is well known.

\begin{theorem}[Kirillov-Kostant-Souriau]\label{Thm:KKS} \cite{K04}
Given any $\xi \in \mathfrak{g}$, 
the adjoint orbit $\mathcal{O}(\xi)$ is a smooth compact connected submanifold in $\mathfrak{g}$ and there exists a natural $G$-invariant (under the adjoint action) symplectic structure on $\mathcal{O}(\xi)$. In other words, there exists a closed non-degenerate $G$-invariant real $2$-form $\omega_{\mathcal{O}(\xi)} \in \Omega^{2}(\mathcal{O}(\xi); \mathbb{R})$ on $\mathcal{O}(\xi)$. More explicitly, $\omega_{\mathcal{O}(\xi)}$ can be constructed in the following way.

For all $\eta \in \mathcal{O}(\xi)$, let $B_{\eta}$ be the antisymmetric bilinear form on $\mathfrak{g}$ defined by
\begin{equation}
B_{\eta}(X, Y) := \left\langle \eta, \left[ X, Y \right] \right\rangle
\end{equation}
for all $X, Y \in \mathfrak{g}$. Then $\omega_{\mathcal{O}(\xi)}$ can be defined by
\begin{equation}
\omega_{\mathcal{O}(\xi)}(\eta)([X,\eta], [Y,\eta]) =     \left\langle \eta, \left[ X, Y \right] \right\rangle 
\end{equation}
for all $X, Y \in \mathfrak{g}$, $\eta \in \mathcal{O}(\xi)$.

Note that for all $\eta \in \mathcal{O}(\xi) \subseteq \mathfrak{g}$, $T_{\eta}\mathcal{O}(\xi) = \left\lbrace [X,\eta] \  : \  X \in \mathfrak{g} \right\rbrace$.

This natural $2$-form $\omega_{\mathcal{O}(\xi)}$ is sometimes referred to as the Kirillov-Kostant-Souriau symplectic form~ on the adjoint orbit $\mathcal{O}(\xi)$.
\end{theorem}

Therefore, an adjoint orbit $\mathcal{O}(\xi)$ becomes a symplectic manifold when it is equipped with its Kirillov-Kostant-Souriau symplectic form $\omega_{\mathcal{O}(\xi)}$. In addition, we have the following:

\begin{proposition}
The adjoint action of $G$ on $(\mathcal{O}(\xi), \omega_{\mathcal{O}(\xi)})$ is a Hamiltonian $G$-action with the moment map given by the inclusion map $\mu_{\mathcal{O}(\xi)}: \mathcal{O}(\xi) \hookrightarrow \mathfrak{g}$. In other words, $\mu_{\mathcal{O}(\xi)}$ is equivariant with respect to the adjoint action of $G$ on $\mathcal{O}(\xi)$ and the adjoint action of $G$ on $\mathfrak{g}$, and for all $X \in \mathfrak{g}$,
\begin{equation}
d \mu_{\mathcal{O}(\xi)}^{X} = \iota_{X^\sharp}\omega_{\mathcal{O}(\xi)}
\end{equation}
where $\mu_{\mathcal{O}(\xi)}^{X}: \mathcal{O}(\xi) \to \mathbb{R}$ is defined by $\mu_{\mathcal{O}(\xi)}^{X}(\eta) = \left\langle \mu_{\mathcal{O}(\xi)}(\eta), X \right\rangle$ for all $\eta \in \mathcal{O}(\xi)$ and $X^\sharp$ is the vector field on $\mathcal{O}(\xi)$ such that for all $\eta \in \mathcal{O}(\xi)$, the tangent vector $X^\sharp(\eta) \in T_{\eta}\mathcal{O}(\xi)$ is
\begin{equation}
\frac{d}{dt} \Bigr|_{t = 0} \left( {\rm Ad} (\exp(tX))\eta \right) \text{.}
\end{equation}
\end{proposition}

Let $\mathcal{O}(\xi_{1}), \cdots, \mathcal{O}(\xi_{N})$ be $N$ adjoint orbits. Then we can form their Cartesian product:
\begin{equation}
\mathcal{M}(\underline{\xi}) := \mathcal{O}(\xi_{1}) \times \cdots \times \mathcal{O}(\xi_{N})
\end{equation}
where
\begin{equation}
\underline{\xi} := (\xi_{1}, \cdots, \xi_{N}) \in \overbrace{\mathfrak{g} \times \cdots \times \mathfrak{g}}^{N} \text{.}
\end{equation}

We assume the following:
\begin{assumption}\label{Assume:G1}
All of $\mathcal{O}(\xi_{1}), \cdots, \mathcal{O}(\xi_{N})$ are diffeomorphic to the homogeneous space $G / T$. This assumption is equivalent to the assumption that all of the stabilizer groups $\Stab_{G}(\xi_{1}), \cdots, \Stab_{G}(\xi_{N})$ are conjugate to the chosen maximal torus $T$. If all of $\xi_{1}, \cdots, \xi_{N}$ are contained in $\mathfrak{t} \subseteq \mathfrak{g}$, then this assumption is saying that 
$$\Stab_{G}(\xi_{1}) = \cdots = \Stab_{G}(\xi_{N}) = T. $$
\end{assumption}

\begin{remark}
	Since every adjoint orbit $\mathcal{O}(\xi)$ can be written as $\mathcal{O}(\xi^{\prime})$ for some $\xi^{\prime} \in \mathfrak{t} \subseteq \mathfrak{g}$, we will always assume that $\underline{\xi} = (\xi_{1}, \cdots, \xi_{N})$ satisfies that $\xi_{j} \in \mathfrak{t} \subseteq \mathfrak{g}$ for all $j$.
\end{remark}

The Cartesian product $\mathcal{M}(\underline{\xi}) = \mathcal{O}(\xi_{1}) \times \cdots \times \mathcal{O}(\xi_{N})$ carries a natural symplectic structure $\omega_{\underline{\xi}}$ defined by:
\begin{equation}
\omega_{\underline{\xi}} := \pr_{1}^{*} \omega_{\mathcal{O}(\xi_{1})} + \cdots + \pr_{N}^{*} \omega_{\mathcal{O}(\xi_{N})}
\end{equation}
where $\pr_{j}: \mathcal{O}(\xi_{1}) \times \cdots \times \mathcal{O}(\xi_{N}) \to \mathcal{O}(\xi_{j})$ is the projection onto the $j$-th component.

Let $G$ act on $\mathcal{M}(\underline{\xi}) = \mathcal{O}(\xi_{1}) \times \cdots \times \mathcal{O}(\xi_{N})$ by the diagonal action $\Delta$:
\begin{equation}
\Delta(g)(\eta_{1}, \cdots, \eta_{N}) := ({\rm Ad}(g)(\eta_{1}), \cdots, {\rm Ad} (g)(\eta_{N}))
\end{equation}
for all $g \in G$, $\eta_{j} \in \mathcal{O}(\xi_{j})$.

We mentioned above  that the symplectic form $\omega_{\underline{\xi}}$ is  invariant under this
 action of $G$.   We also have the following:
\begin{proposition}
The diagonal action $\Delta$ of $G$ on $(\mathcal{M}(\underline{\xi}), \omega_{\underline{\xi}})$ is a Hamiltonian $G$-action with the moment map $\mu_{\underline{\xi}}: \mathcal{M}(\underline{\xi}) \to \mathfrak{g}$ being:
\begin{equation}
\mu_{\underline{\xi}}(\underline{\eta}) = \sum_{j = 1}^{N} \eta_{j}
\end{equation}
for all $\underline{\eta} := (\eta_{1}, \cdots, \eta_{N}) \in \mathcal{M}(\underline{\xi})$.
\end{proposition}

We assume that:
\begin{assumption}\label{Assume:G2}
$0 \in \mathfrak{g}$ is a regular value for $\mu_{\underline{\xi}}: \mathcal{M}(\underline{\xi}) \to \mathfrak{g}$ and $\mu_{\underline{\xi}}^{-1}(0) \neq \emptyset$.
\end{assumption}

\begin{remark}
By Sard's theorem, the set
\begin{equation}
\mathcal{A} := \left\lbrace \underline{\xi} \in \overbrace{\mathfrak{t} \times \cdots \times \mathfrak{t}}^{N} \  : \  \text{Assumptions \ref{Assume:G1}, \ref{Assume:G2} hold} \right\rbrace
\end{equation}
has nonempty interior in $\mathfrak{t} \times \cdots \times \mathfrak{t}$.
\end{remark}

Then, the level set $$\mathcal{M}_{0}(\underline{\xi}) := \mu_{\underline{\xi}}^{-1}(0)$$ is a closed, thus compact, submanifold of $\mathcal{M}(\underline{\xi})$ and the diagonal action $\Delta$ of $G$ restricts to an action on $\mathcal{M}_{0}(\underline{\xi})$. Therefore, we can form the quotient space (or symplectic reduction)  with respect to this action of $G$ on $\mathcal{M}_{0}(\underline{\xi})$:
\begin{equation}
\mathcal{M}_{\red}(\underline{\xi}) := \mathcal{M}_{0}(\underline{\xi}) / G \text{.}
\end{equation}
The quotient space is also compact.

If the $G$-action on $\mathcal{M}_{0}(\underline{\xi})$ is free and proper (in our situation, properness is automatically satisfied), then the quotient space $\mathcal{M}_{\red}(\underline{\xi}) = \mathcal{M}_{0}(\underline{\xi}) / G$ is a smooth manifold. However, in our situation, the $G$-action on $\mathcal{M}_{0}(\underline{\xi})$ is in general not free. Hence in general 
it follows from the treatment in \cite{GS} that the quotient space is  an orbifold \cite{HM} rather than a smooth manifold. To avoid this complication, we will make the following assumption.
\begin{assumption}\label{Assume:G3}
The quotient space $\mathcal{M}_{\red}(\underline{\xi}) = \mathcal{M}_{0}(\underline{\xi}) / G$ is a smooth compact manifold.
\end{assumption}

Assumption \ref{Assume:G3} is satisfied provided the stabilizer of the action of $G$ at 
all points in  $\mathcal{M}_{0}(\underline{\xi})$ is the identity. 

\begin{remark}
The above assumption will put further restrictions on which $\underline{\xi} \in \mathfrak{t} \times \cdots \times \mathfrak{t}$ we can choose as initial data. Thus we only choose initial data from the following set:
\begin{equation}
\mathcal{A}^{\prime} := \left\lbrace \underline{\xi} \in \overbrace{\mathfrak{t} \times \cdots \times \mathfrak{t}}^{N} \  : \  \text{Assumptions \ref{Assume:G1}, \ref{Assume:G2}, and \ref{Assume:G3} hold} \right\rbrace
\end{equation}

%
Notice that since the elements in the center of $G$ always act trivially on $\mathcal{M}(\underline{\xi})$ and $\mathcal{M}_{0}(\underline{\xi})$, Assumption \ref{Assume:G3} is valid if $PG = G / Z(G)$ acts freely on $\mathcal{M}_{0}(\underline{\xi})$. This happens for $G = \SU(n)$ if all the coadjoint orbits $\mathcal{O}(\xi_{i})$ are generic.
\end{remark}

Then, we have the following well known theorem:
\begin{theorem}[Marsden-Weinstein]

The smooth compact manifold 
$$\mathcal{M}_{\red}(\underline{\xi}) = \mathcal{M}_{0}(\underline{\xi}) / G$$ carries a unique symplectic structure $\omega_{\red}(\underline{\xi})$ such that
\begin{equation}
i^{*} \omega_{\underline{\xi}} = \pi^{*} \omega_{\red}(\underline{\xi})
\end{equation}
where $i: \mathcal{M}_{0}(\underline{\xi}) \hookrightarrow \mathcal{M}(\underline{\xi})$ is the inclusion map and $\pi: \mathcal{M}_{0}(\underline{\xi}) \to \mathcal{M}_{\red}(\underline{\xi})$ is the associated projection map.
\end{theorem}

\begin{definition}
We call this compact symplectic manifold $$(\mathcal{M}_{\red}(\underline{\xi}), \omega_{\red}(\underline{\xi}))$$
 an \emph{{$N$-fold reduced product}}.
\end{definition}

%
%
%

\begin{remark}
The dimension of an $N$-fold reduced product is
\begin{equation}
N (\dim G - \dim T) - 2 \dim G = (N-2) \dim G - N \dim T 
\end{equation}
when all  orbits are generic. 
In the case $G = \SU(3)$ and $N = 3$ , this is $\dim G - 3 \dim T = 8 - 6 = 2$. 
These reduced products are diffeomorphic to the 2-sphere \cite{TRP1}.
\end{remark}


%

\begin{remark}
If the initial point $\underline{\xi}$ is clear from the context, we will suppress the inclusion of the point $\underline{\xi}$ in our notations and write, for example, $\mathcal{M}, \mathcal{M}_{0}, \mathcal{M}_{\red}$ instead of $\mathcal{M}(\underline{\xi}), \mathcal{M}_{0}(\underline{\xi}), \mathcal{M}_{\red}(\underline{\xi})$, respectively. Similarly, this is done for the notations of the symplectic structures and so on.
\end{remark}

\section{Intersection pairings of $N$-fold reduced products}

\subsection{Introduction}
In our previous paper \cite{JJ18a}, we investigated the symplectic volume of $N$-fold reduced products and derived the following formula for all generic $N$-fold reduced products:

\begin{theorem} \label{th9} 
In the notation introdued earlier, and under the hypotheses imposed in the previous section, we  have
	\begin{equation} \label{eqsixteen}
	\int_{\mathcal{M}_{\red}} e^{i \omega_{\red}} = \sum_{\underline{w} \in W^{N}} \sgn(\underline{w}) \int_{X \in \liet} \frac{e^{i \left\langle \mu_{T}(\underline{w} \cdot \underline{\xi}), X \right\rangle}}{\varpi^{N - 2}(X)} \,  dX
	\end{equation}
	where $\mu_{T}: \mathcal{M} \to \liet$ is the moment map for the $T$-action on $\mathcal{M}$, $\underline{\xi} = (\xi_{1}, \cdots, \xi_{N}) \in (\liet)^{N}$ is generic, $\underline{w} = (w_{1}, \cdots, w_{N}) \in W^{N}$ and
	\begin{equation}
	\varpi(X) = \prod_{\gamma} \left\langle \gamma, X \right\rangle
	\end{equation}
	where $\gamma$ runs over all the positive roots of $G$.
\end{theorem}

\subsection{Equivariant cohomology and the Cartan model}
The main tool we used to prove Theorem \ref{th9} is the Atiyah-Bott-Berline-Vergne localization
formula.  (See \cite{JK95}.)
We make use of the Cartan model for 
equivariant cohomology (see for example \cite{Meinrenken}).
In this model, an equivariant differential form 
is represented by a linear combination of differential forms $\alpha_j$
with polynomial dependence on a parameter $$X \in \liet.$$ 
We assume $\alpha_j $ has degree $j$ in $X$.
The grading is the sum of the  differential form grading and two times  the degree as a polynomial in $X$.
The differential is $$d_X =d - \iota_{X^{\sharp}} $$
where $\iota$ denotes interior product.
Recall that  $X^\sharp$ is the fundamental vector field generated by the 
action of $X$.
For example, the extension of the symplectic form to an equivariantly closed form is
$$\bar{\omega}(X) = \omega + \mu_{X} $$ where
$\mu_X$ is the moment map associated to $X$ (in other words the function whose Hamiltonian 
vector field is $X^\sharp$).

An equivariant $m$-form  $\alpha $ in the Cartan model is 
a sum of terms $\alpha_j$ for $2j \le m$, where the degree of $\alpha_j$ as a differential form
is $m-2j$. 
If the differential form degree is $0$, then $j = m/2$ where $m$ is the (real) dimension 
of the manifold. 

The restriction of $\alpha$ to a fixed point of the $T$ action is 
$\alpha_{m/2} $ (the term of degree $0$ as a differential form).
If the form $\alpha $  is equivariantly closed, 
it follows that
$$ d \alpha_j = \iota_{X^\sharp} \alpha_{j-1}$$ for all $j$.


Let $M$ be a Hamiltonian $G$-manifold.
The Kirwan map, which we shall denote by $\kappa$,  is a map 
 from $H^*_G(M) $ to $H^*_G(M_0)$, where 
$M_0$ is defined as the zero level set of the moment map on $M$.   It is  the restriction map
to a level set of the moment  map.  If $0$ is a 
regular value of the moment map, then $H^*_G(M_0) 
\cong H^*(M_0/G)$ .  When   $0$ is a regular value of the moment map,
Kirwan proved that the map $\kappa$  is surjective \cite{Kthesis}.


\subsection{Cohomology of orbits}
For an adjoint orbit homeomorphic to $G/T$, we see (for example from \cite{F}, Chap. 10.2 Proposition 3)
that  the cohomology is generated multiplicatively 
 by the first Chern classes of line bundles $L_\beta$ over the orbit,
where 
\begin{equation}
L_\beta = G \times_{T, \beta} \CC
\end{equation}
where we write the orbit as $G/T$ and 
the equivalence relation is $$(g,z) \sim (gt, \beta(t)^{-1} z) $$
for $g \in G$, $t \in T$, $z \in \CC$  and 
for a weight $\beta \in {\rm Hom}(T, U(1)). $
For example, for $G= SU(n)$,  the collection of $\beta$ comprising  the simple roots of $G$
gives rise to a basis for the cohomology of $G/T$.  For $G = SU(n)$, 
a proof of this result can be found  in  Fulton's book \cite{F}
(Chapter 10.2,  Proposition 3). For general Lie groups
this is  Theorem 5 in Section 4 in  the article  by  Tu \cite{Tu}.

We can write each weight $\beta$ as 
$$ \beta (\exp X) = \exp (2 \pi  B(X))$$
for a linear map $B: \liet \to \RR$ which sends the integer  lattice (the kernel of the exponential map)  to $\ZZ$. 
Here we have used the exponential map $\exp: \liet \to T$. 
The equivariant first Chern class of the line bundle $L_{\beta}$ is denoted
$${c_1}^{\rm eq} (L_b).$$
Its restriction to an isolated fixed point $F$
is $$ {c_1}^{\rm eq} (L_b)|_F =  c_1(L_\beta)|_F  + B(X). $$

The restriction of this equivariant first Chern class to a component $F$ of the fixed point set is $B(X)$.
By naturality,
we have that 
\begin{equation}
\pi_j^* \bigl ( c_1 (L_j) \bigr )  = c_1 \bigl (\pi_j^* L_j\bigr) 
\end{equation}
where $$ \pi_j:  {\mathcal{O}}_{\xi_1} \times \dots \times {\mathcal{O}}_{\xi_N} 
\to {\mathcal{O}}_{\xi_j}$$
is projection on the $j$-th orbit, and $L_j$ is a line bundle over 
${\mathcal{O}}_{\xi_j}$.

\subsection{Localization}
The Atiyah-Bott-Berline-Vergne localization formula 
leads to the following (see \cite{JK95}, Theorem 8.1):
$$ \int_{M_{\rm red}} \kappa(\alpha) = {\rm Res} \sum_F \alpha_{m/2} (X) \frac{e^{i \mu_X(F)} } {e_F(X) }.
$$

In the case when $M$ is the product of $N$ adjoint orbits
when $$\alpha =  \exp (i\bar{\omega}) $$ is the equivariant extension of the
symplectic volume form, and 
$$\kappa (\alpha) =  e^{i \omega_{\rm red} }$$ 
is the symplectic volume form on $M_{\rm red}$.   Theorem \ref{th9} may be expressed as follows. 

\begin{equation}  \label{eqsixteenb}
\int_{M_{\rm red}}  \kappa(\alpha)  = {\rm Res} \sum_{w \in W} {e^{i (w \lambda,X) }}
\sgn(w) \frac{1}{ \bigl ({\varpi}(X) \bigr )^{N-2}  }. \end{equation}
Equation (\ref{eqsixteenb}) is the meaning
of 
 the integral over $\liet$  in equation (\ref{eqsixteen}) 
 whose definition is given in \cite{GLS}
and elaborated in \cite{JK95}.
The symbol ${\rm Res} $ (the residue) is defined in \cite{JK95}, Theorem 8.1. See also 
\cite{JK3}, Proposition 3.2. 
The residue  has several
equivalent definitions (as outlined in \cite{JK3}). One of these
definitions characterizes the residue  as  an iteration of    one-variable residues.

\begin{remark}
One feature that is special to our situation (Cartesian products of adjoint 
orbits) is that all the equivariant Euler classes are the same,
except for the sign (which is $\sgn(\underline{w})$, the product of the  signatures of 
the 
permutations). Up to  sign, the equivariant Euler class is a power $\varpi(X)^N$  of 
$\varpi(X)$ where $\varpi$ is the product of positive roots.
\end{remark}

In the above notation, we have the following generalization of Theorem \ref{th9}:

\begin{theorem}\label{theorem:inteq}
Let $\mathcal{M}$ be as above, and
let $\zeta$ be a $G$-equivariant cohomology class on $\mathcal{M}$.  Let $\kappa: H^*_G(\mathcal{M}) \to 
H^*(\mathcal{M}_{\red}) $ be the
Kirwan map. We have
\begin{equation}
\int_{\mathcal{M}_{\red}} e^{i \omega_{\red}}\kappa(\zeta)  
= \sum_{\underline{w} \in W^{N}} \sgn(\underline{w}) \int_{X \in \liet} \frac{e^{i \left\langle \mu_{T}(\underline{w} \cdot \underline{\xi}), X \right\rangle}
\zeta(X)|_{\underline{w} \cdot \underline{\xi}}}{\varpi(X)^{N - 2}} \,  dX.
\end{equation}

\begin{equation} \label{eqnineteenb}
= {\rm Res} \sum_{\underline{w} \in W^N} {e^{i (\underline{w} \cdot  \underline{\xi},X) }}
\sgn(\underline{w}) \frac{\zeta(X)|_{\underline{w} \cdot \underline{\xi} } }{ \bigl ({\varpi}(X) \bigr )^{N-2}  }. \end{equation}
Here $\zeta(X)$ is a product of powers of  a collection of 
equivariant first Chern classes $\left ( c_1^{\rm eq} (L_{{\beta_\ell} }(X)\right )^{n_\ell}  $ where the index $\ell $ runs from $1$ to $ N$ if we are considering the reduced space of the 
product of $N$ orbits and $n_\ell$ is a nonnegative integer, and
the weight of the $\ell$-th line bundle is $\beta_\ell$ with associated 
linear map $B_\ell$.
The restriction  of $\zeta$ 
to the fixed point set of the 
$T$ action is 
$$\prod_{\ell}  \left ( B_{\ell} (X) \right )^{n_\ell}  . $$



\end{theorem}

\begin{remark}
Theorem \ref{theorem:inteq} describes  all intersection pairings between cohomology classes 
of $\mathcal{M}_{\red}$. 
\end{remark}

\end{document}